\magnification=\magstep1 \voffset0.4in
\vsize=8.9truein
\hsize=138truemm
\overfullrule=0pt      \input epsf  
\nopagenumbers 
\headline={\ifnum\pageno=1 \hfil\else\ifodd\pageno\rightheadline \else\leftheadline\fi\fi}
\def\rightheadline{\tenrm\hfil CFT and mapping class groups\hfil\folio}
\def\leftheadline{\tenrm\folio\hfil  CFT and mapping class groups\hfil}
\voffset=2\baselineskip
 \def\leaderfill{\leaders\hbox to 1em{\hss.\hss}\hfill} 
 \def\la{\lambda}   
       \def\i{{\rm i}}

   \def\H{{\Bbb H}}
\font\huge=cmr10 scaled \magstep2

\font\smcap=cmcsc10        

\input amssym.def
\def\Z{{\Bbb Z}} \def\R{{\Bbb R}} \def\Q{{\Bbb Q}}  
\def\C{{\Bbb C}}

\def\buildrel#1\over#2{{\mathrel{\mathop{\kern0pt #2}\limits^{#1}}}}

\centerline{{\bf \huge Conformal field theory and mapping class groups}}
\bigskip

\centerline{{Terry Gannon}}\smallskip

\centerline{{Mathematics Dept, U of Alberta, Canada}}\bigskip

{\narrower\noindent{\bf Abstract.} Rational conformal field theories
produce a tower of finite-dimensional representations of surface 
mapping class groups, acting on the conformal blocks of the theory.
We review this formalism. We show that many recent mathematical developments 
can be fit into the first 2 floors of this tower. We also review what is known
in higher genus.

}\bigskip

\noindent{{\bf 1. Introduction.}}

Conformal field theory (CFT) is an exceptionally symmetric
quantum field theory. It takes classical mathematical structures, and 
`loops' or `complexifies' them, to produce 
infinite-dimensional structures such as affine algebras or vertex operator
algebras. Its importance to math -- which is considerable -- is that the
resulting  structures tend to straddle several areas, such as geometry,
 algebra, number theory, functional analysis, ....

One of its most beautiful but least appreciated accomplishments is the
reorganisation of several recent mathematical developments, such as
Monstrous Moonshine, Jones' knot invariant,
the modularity of affine Kac-Moody characters, and 
 braid group monodromy of the KZ equation, 
into the first two floors of an infinite tower. This paper describes the
resulting picture. For more details see the book  [18] and references therein.

\bigskip\noindent{{\bf 2. Mapping class groups.}}

Up to homeomorphism, a (connected oriented real) surface is completely
characterised by its genus $g$ and number $n$ of boundary components (punctures). In fact,
 there's only one way (up to equivalence) to give it a
real-differential structure.
But a real surface can also be a complex curve -- it can usually be
given a complex-differential (equivalently, a conformal) structure in
infinitely many different ways.

For example,  the torus $\R^2/\Z^2$ can be given a complex
structure by replacing $\R^2$ with $\C$ and $\Z^2$ by 
$\Z+\tau\Z$ for any $\tau\in\C$
with nonzero imaginary part. In fact, any torus is 
conformally equivalent to one of the form $\C/(\Z+\tau\Z)=:T_\tau$, where
$\tau$ lies in the {upper-half plane}
$\H:=\{x+\i y\in\C\,|\,y>0\}$.
Moreover, the tori $T_\tau$ and $T_{\tau'}$ are themselves conformally equivalent, iff
$\tau'={a\tau+b\over c\tau+d}$
for a matrix $\left(\matrix{a&b\cr c&d}\right)\in {\rm SL}_2(\Z)$.

The set of possible complex structures on the torus forms the
{\it moduli space} ${\frak M}_{1,0}$, so labelled because the torus has
genus 1 and 0 punctures. This moduli space can be identified with
the orbifold $\H/{\rm SL}_2(\Z)$:
we call SL$_2(\Z)=:\Gamma_{1,0}$ its {\it mapping class
group}, and  $\H=:{\frak T}_{1,0}$ its {\it Teichm\"uller
space}. Similarly, $\Gamma_{g,n}$,
 ${\frak T}_{g,n}$, and ${\frak M}_{g,n}={\frak T}_{g,n}/\Gamma_{g,n}$ can be
defined for any other $g,n\ge 0$. In particular,
the Teichm\"uller space ${\frak T}_{g,n}$ (a simply connected
complex manifold) accounts for `continuous'
conformal equivalences, while the mapping class group $\Gamma_{g,n}$
(almost always an infinite discrete group) contains the left-over `discontinuous'
ones. The moduli spaces usually have conical singularities, corresponding
to surfaces with extra symmetries; taking into account these stabilisers,
$\Gamma_{g,n}$ will be the (orbifold) fundamental group of ${\cal M}_{g,n}$.

For example,
${\frak M}_{1,0}$ is a sphere with a puncture (corresponding to the cusp
$\Bbb{Q}\cup \i\infty$), and 
conical singularities at $\tau=\i$ and $e^{2\pi\i/3}$.
Because $\C/(\Z+\tau\Z)$ can also be interpreted as a torus with a special point, namely the
additive identity 0, we also have ${\frak T}_{1,1}=\H$ and
$\Gamma_{1,1}={\rm SL}_2(\Z)$.

The surfaces relevant to our story possess additional structure.
Let $\Sigma$ be a compact genus-$g$ surface
with $n$ marked points $p_i\in\Sigma$. About
each point $p_i$ choose a local coordinate $z_i$, vanishing at $p_i$ --
this identifies a neighbourhood
of $p_i$ with a neighbourhood of $0\in\C$. We call $(\Sigma,\{p_i\},\{z_i\})$ an {\it enhanced
surface}{} of type $(g,n)$. The resulting moduli space $\widehat{{\frak M}}_{g,n}$
is infinite-dimensional, but its mapping class group $\widehat{\Gamma}_{g,n}$ is an extension
of $\Gamma_{g,n}$ by $\Z^n$.
For example, $\widehat{\Gamma}_{1,1}$ is the braid group ${\cal B}_3$.

As we will see below, a 
rational conformal field theory gives finite-dimensional representations
of each $\widehat{\Gamma}_{g,n}$ -- merely projective for $n=0$, but it seems
a true one for $n\ge 1$ (though a proof of trueness is to my knowledge only available
for $g\le 1$). Enhanced surfaces are important because they have canonical sewings.
Nevertheless it is common to restrict instead to the projective
representations of $\Gamma_{g,n}$, and pay at most lipservice to the coordinates
$z_i$.

\bigskip\noindent{{\bf 3. Conformal field theory.}}

This section introduces the correlation functions and chiral blocks of conformal
field theory.

A {conformal field theory} (CFT) is a quantum field theory, usually on a two-dimensional space-time
$\Sigma$, whose symmetries include the conformal transformations 
(so conformally equivalent space-times are identified). 
 We restrict to compact orientable $\Sigma$. The same
CFT lives simultaneously on all such $\Sigma$.  See e.g. [9,16,19,27], 
and Chapter 4 of [18] for reviews.

Two dimensions are special for CFT because the local
conformal maps, which form the Lie algebra ${\frak{so}}_{n+1,1}
(\R)$ in $\R^n$ for $n>2$, becomes infinite-dimensional in $\R^2$
(thanks to their identification with (anti-)holomorphic maps). 
 The conformal algebra in two dimensions consists of two commuting 
copies of the Witt algebra $\frak{Witt}$ (one for the
holomorphic maps, and the other for anti-holomorphic ones).
$\frak{Witt}$ is the infinite-dimensional Lie algebra of vector fields
on $S^1$, and has a basis $\ell_{n}$, $n\in\Z$, satisfying
$$[\ell_m,\ell_n]=(m-n)\ell_{m+n}\ .\eqno(3.1)$$
Its unique nontrivial central extension is the Virasoro algebra $\frak{Vir}$,
with basis $L_n,C$ satisfying $[L_n,C]=0$ and
$$[L_m,L_n]=(m-n)L_{m+n}+\delta_{n,-m}{m\,(m^2-1)\over 12}\,C\ .
\eqno(3.2)$$

Basic data in the CFT are the quantum fields $\varphi(z)$,
called {\it vertex operators}, living on space-time $\Sigma=S^2=\C\cup\{\infty
\}$ and centred at $z=0$. Being quantum
fields, these $\varphi$ are 
`operator-valued distributions' on $\Sigma$, acting on the space ${\cal H}$ of 
states for $\Sigma$. The most important vertex operators 
 are the {\it stress-energy tensors}
$T(z),\overline{T}(z)$, 
which are the conserved currents of the conformal symmetry,
as promised by Noether's Theorem; the corresponding conserved charges are
operators $L_n,\overline{L}_m$ defining a $\frak{Vir}$-representation, with central term $C$ 
given by scalars $cI,\overline{c}I$ called the {\it central charges}.

In a typical quantum field theory, a theoretical physicist
makes contact with experiment by computing transition amplitudes 
$\langle{\rm out}|{\rm in}\rangle$ between in-coming and out-going states,
 given mathematically by a Hermitian product $(|{\rm
out}\rangle,S|{\rm in}\rangle)$ in the Hilbert space ${\cal H}$ of states,
for some operator $S$ called the scattering matrix. In practise these can 
 only be calculated in infinite time ($t\rightarrow\pm \infty$) 
limits. The typical way (`LSZ reduction formulae') to express these
asymptotic amplitudes is via artifacts sometimes called {\it correlation
functions}. The theory is regarded as solved if all correlation
functions can be computed. We are interested in the correlation functions
$$\left\langle \varphi_{1}(z_1)\,\varphi_{2}(z_2)\cdots\varphi_{n}(z_n)
\right\rangle_{\Sigma;p_1,\ldots,p_n}\eqno(3.3)$$
 of CFT, for any vertex operators $\varphi_i$ and
any enhanced surface $(\Sigma,\{p_i\},z_i)$ (so
 $\varphi_{i}(z_i)$ is `centred' at $p_i\in\Sigma)$. 
The remainder of this section explains how physicists think of these
CFT correlation functions. Their intuition is provided by string theory.

In a typical quantum theory, correlation functions
are calculated perturbatively by Taylor-expanding in some coupling constant.
For this purpose,
Feynman's path integral formulation --  the quantisation of Hamilton's 
action principle
in classical mechanics -- is convenient.
Each term in this perturbation series is computed
separately using Feynman diagrams and rules. 
 The Feynman diagrams of quantum field theory
are graphs, with a different kind of edge for each species of particle, and
a different kind of vertex for every term in the interaction part of the Lagrangian;
Feynman's rules describe how to go from these diagrams
 to certain integral expressions and hence to the individual terms
 in the Taylor series expansion of the given correlation function.
Feynman diagrams are combinatorial artifacts describing `virtual' (non-real)
processes; topologically equivalent ones are identified,
and in practise only the simplest  are ever considered.

Applying this perturbation formalism to string theory 
recovers CFT. Consider for convenience closed strings. Then
CFT lives on the world-sheet
$\Sigma$ (string theory's Feynman diagrams) traced by the strings as they
virtually evolve, colliding and separating, through time:
string amplitudes (in e.g.\  26-dimensional space-time) 
can be expressed as correlation functions of a (point-particle) CFT
in two dimensions. The boundaries
of these world-sheets are the in-coming and out-going strings; 
the world-sheets for asymptotic amplitudes have semi-infinite end-tubes and
 can be conformally mapped to compact surfaces with punctures $p_i$
(one for every external string). The corresponding Feynman integral
is over moduli space $\widehat{{\cal M}}_{g,n}$. The data
of those external strings are stored in the appropriate vertex operator
attached to that point $p_i$. 
The Witt algebra arises here as infinitesimal reparametrisations
of the string (a circle).

Everything in CFT comes in
 a combination of strictly holomorphic, and strictly anti-holomorphic,
 quantities. Here,
`holomorphic' is in terms of space-time $\Sigma$ (which locally looks
like $\C$), or the appropriate
moduli space (which usually locally looks like $\C^\infty$). These holomorphic and
anti-holomorphic building blocks are called {\it chiral}. A CFT is studied by
first analysing its chiral parts, and then determining explicitly how they
piece together to form the physical(=bi-chiral) quantities.
Almost all attention by
mathematicians has focused on the chiral (as opposed to physical) data.
In string theory, this holomorphic/anti-holomorphic alternative corresponds to
classical ripples travelling clockwise/anti-clockwise around the string.

A typical vertex operator $\varphi(z)$ depends neither holomorphically nor
anti-holomorphically on $z$.
Let ${\cal V}$ consist of all the holomorphic, and
$\overline{\cal V}$ the anti-holomorphic, ones.
For example, the {stress-energy tensor} $T(z)$ and all of its derivatives
lie in ${\cal V}$.
In the very simplest CFTs, called the minimal models, ${\cal V}$
 consists only of $T(z)$ and its derivatives.

These {\it chiral algebras} ${\cal V},\overline{\cal V}$ have a
rich mathematical structure, with a `multiplication' coming from the so-called
operator product expansion, and are examples of {\it vertex operator algebras}
(see  [28], or Chapter 5 of  [18]).
${\cal V}$ and $\overline{\cal V}$  mutually commute and
the full symmetry `algebra' of the CFT can be identified with 
${\cal V}\oplus\overline{\cal V}$.
Since quantum fields act on state-space ${\cal H}$, it
carries a representation of ${\cal V}\oplus\overline{
{\cal V}}$ and decomposes into a direct integral  of 
irreducible  ${\cal V}\oplus\overline{{\cal V}}$-modules. 
A {\it  rational conformal field theory} (RCFT) is one whose 
state-space ${\cal H}$ decomposes in fact into a {\it finite direct sum}
$${\cal H}=\oplus_{M\in\Phi,\overline{N}\in\overline{\Phi}}\,
{\cal Z}_{M,\overline{N}}\,M\otimes \overline{N}\ ,\eqno(3.4)$$
where $\Phi$ and $\overline{\Phi}$ denote the (finite) sets of
irreducible ${\cal V}$- and $\overline{\cal V}$-modules,
and the ${\cal Z}_{M,\overline{N}}\ge 0$ are multiplicities.
The RCFT  are especially symmetric and well-defined quantum field theories 
and are the CFTs we're interested in.
 The name `rational' arises 
because their central charges $c,\overline{c}$ lie in $\Q$.

The correlation functions (3.3) can be expressed in terms of purely
chiral
quantities called {\it conformal} or {\it chiral blocks}, denoted
$$\left\langle {\cal I}_1(v_1,z_1)\,{\cal I}_2(v_2,z_2)\cdots
{\cal I}_n(v_n,z_n)\right\rangle_{(\Sigma;p_1,\ldots,p_n;M^1,\ldots,M^n)}\ .
\eqno(3.5)$$
As usual, $(\Sigma,\{p_i\},\{z_i\})$ is an enhanced surface, and
to each $p_i$ we assign a ${\cal V}$-module $M^i\in \Phi$ and a 
 state $v_i\in M^i$. The holomorphic field ${\cal I}_i(v_i,z_i)$ centred at 
$p_i$ is an operator-valued distribution called an {\it intertwiner}
sending (`intertwining') one ${\cal V}$-module (say $M\in
\Phi$) to another (say  $N\in\Phi$). For fixed $M,N,M^i$,
the dimension of the space of intertwiners
is called the {\it fusion coefficient} ${\cal N}_{M^i,M}^N$, and is given by
Verlinde's formula 
$${\cal N}_{M^i,M}^N=\sum_{P\in \Phi}{S_{M^iP}S_{MP}\overline{S_{NP}}
\over  S_{{\cal V}P}}\ ,\eqno(3.6)$$
where the matrix $S$ (no relation to the scattering matrix) is defined in
(5.1a) below. For example, ${\cal V}$ is always a module
for itself and ${\cal N}_{{\cal V},{\cal V}}^{{\cal V}}=1$; the unique
(up to scaling) intertwiner ${\cal I}(v,z)$ bijectively associates states
$v$ with vertex operators $\varphi(z)={\cal I}(v,z)$ (the so-called 
`state-field correspondence'). Thus intertwiners generalise vertex operators 
$\varphi\in{\cal V}$.

To solve a given RCFT, the strategy then is to:

\smallskip \item{(a)} construct all possible chiral
blocks (3.5); and

\smallskip\item{(b)} construct the correlation
functions (3.3) from those chiral blocks.

\bigskip\noindent{\bf 4. The chiral blocks of RCFT.}

For a fixed $(g,n;M^1,\ldots,M^n)$, an RCFT assigns
a finite-dimensional space ${\frak F}^{(g,n)}_{(M^i)}$  
of chiral blocks (3.5). Chiral blocks are important to
RCFT because finite combinations of them are the correlation functions, and 
knowing the latter is equivalent to solving the theory.

Each chiral block depends multi-linearly
on the states $v_i\in M^i$, and holomorphically on the $z_i$, provided branch-cuts
in $\Sigma$ between $p_i$ are made; locally, it can be regarded as a 
holomorphic function on $\widehat{{\cal M}}_{g,n}$. The dimension of this
space is given by Verlinde's formula
$${\rm dim}\,{\frak F}^{(g,n)}_{(M^i)}=\sum_{P\in\Phi}
{S_{M^1 P}\over S_{{\cal V} P}}\cdots {S_{M^n P}\over S_{{\cal V} P}}
S_{{\cal V}P}^{2(1-g)}\ ,\eqno(4.1)$$
a generalisation of (3.6).

Moore and Seiberg [30] -- see also [2] -- isolated 
the data (finite-dimensional vector spaces and linear transformations)
defining each chiral half of RCFT, and provided a complete set of relations they satisfy. 
Huang is pursuing the explicit construction
 for all sufficiently nice chiral algebras ${\cal V}$ (see e.g.\ [21]).

A basis for ${\frak F}^{(g,n)}_{(M^i)}$ is found by performing the following 
Feynman rules (called `conformal bootstrap'). Fix a surface $\Sigma$
  of type $(g,n)$. The space ${\frak F}^{(0,3)}_{(M,N,P)}$ consists of
intertwiners; arbitrarily fix bases for all those spaces. Now,
dissect $\Sigma$ into pairs-of-pants, as in Figure 1;
assign  a dummy label $N_j\in\Phi$ to each internal cut in the  dissection;
to each vertex in your dissection, 
 choose an intertwining operator from the
basis of the appropriate space ${\frak F}^{(0,3)}$;
`evaluate' the corresponding chiral block  -- e.g.\ for
each cut, a trace is taken of the product of intertwiners.
Repeating, by running through all possible values of the dummy labels,
the result is a basis of chiral blocks.

\medskip\epsfysize=1in\centerline{ \epsffile{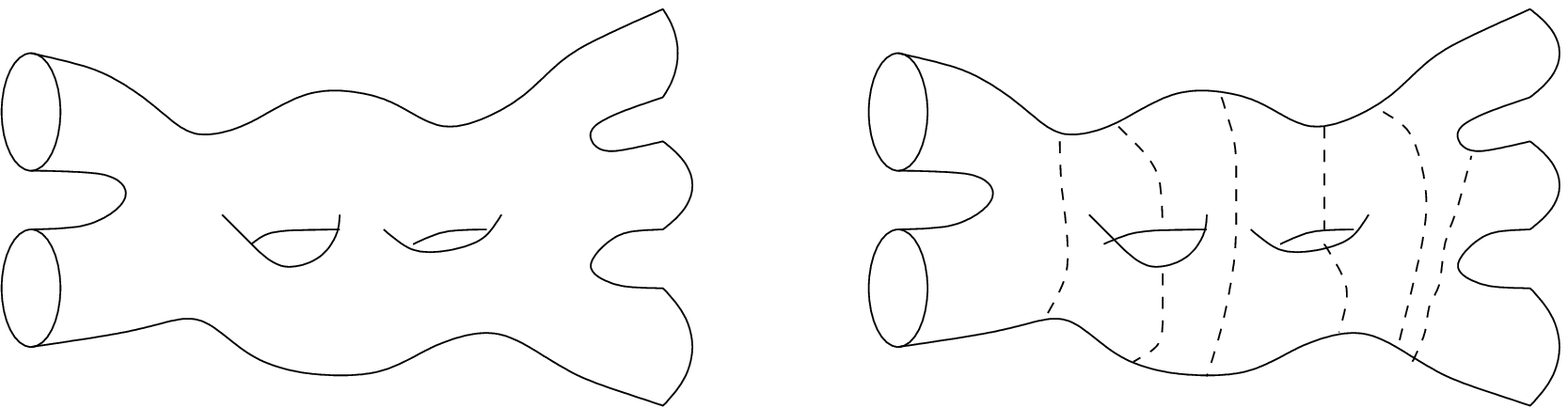}}\medskip
\centerline{{\smcap Figure 1.} Dissecting a surface into pairs-of-pants}\medskip

For an important example, let $\Sigma$ be the torus $\C/(\Z+\tau\Z)$ with one
puncture $p$ (say at 0),
assigned the module $M^1={\cal V}$ 
and state $v_1=|0\rangle$ (the vacuum, the state of lowest energy). One cut 
suffices to unfold it into a sphere with 3 punctures, assigned ${\cal V}$-modules
${\cal V},M,M\in\Phi$ ($M$ is the dummy label). The fusion coefficient 
${\cal N}_{{\cal V},M}^M$ always  equals 1, and so for each
$M\in\Phi$ there is a unique intertwiner, say ${\cal I}_1^{(M)}$.
Hence dim$\,{\frak F}^{(1,1)}_{{\cal V}}=\|\Phi\|$.
These Feynman rules yield the {chiral block} 
$$\chi_M(\tau):={\rm tr}_M e^{2\pi\i\tau\,(L_0-c/24)}\ ,\eqno(4.2)$$
where $c$ is the central charge and
the Virasoro element $L_0$ corresponds to energy (the trace comes from the 
dissection). These span ${\frak F}^{(1,0)}$.

One of the simplest RCFT is the Ising model, a minimal model. It has central
charge $c=\overline{c}=0.5$, and its
chiral algebra has 3 irreducible modules, which we'll label $\Phi=
\overline{\Phi}=\{{\cal V},\epsilon,\sigma\}$. Its toroidal chiral blocks
(4.2) are
$$\eqalignno{\chi_{{\cal V}}(\tau)=&\,q^{-1/48}\,(1+q^2+q^3+2q^4+2q^5+3q^6+3q^7+\cdots)
\ ,&\cr
\chi_\epsilon(\tau)=&\,q^{23/48}\,(1+q+q^2+q^3+2q^4+2q^5+3q^6+3q^7+\cdots)
\ ,&(4.3)\cr
\chi_\sigma(\tau)=&\,q^{1/24}\,(1+q+q^2+2q^3+2q^4+3q^5+4q^6+5q^7+\cdots)\ ,
&}$$
where $q=e^{2\pi\i\tau}$. 

Other important RCFT are the {\it Wess--Zumino--Witten}{} (WZW{})
models. These correspond to strings living on a compact Lie group $G$. 
The chiral algebra ${\cal V}$ is closely
related to the affine Kac--Moody algebra ${\frak{g}}^{(1)}$ (see
[23]), where
${\frak{g}}$ is the Lie algebra of $G$ (${\frak{g}}^{(1)}$
is the nontrivial central extension of the loop algebra ${\frak{g}}
\otimes\C[z^{\pm 1}]$).
Its modules $M\in\Phi$ can be identified with the integrable highest-weight 
modules $L(\lambda)$ at a level $k$ determined by the central charge $c$.
The chiral blocks $\chi_M(\tau)$ for the WZW models are a specialisation of
the corresponding affine algebra ${\frak{g}}^{(1)}$-character $\chi_\lambda(h)$, 
and for this reason $\chi_M(\tau)$ in any RCFT is called the character of the 
${\cal V}$-module $M$.  We'll return to the WZW and Ising models shortly.

Each dissection produces a basis for the space ${\frak F}^{(g,n)}_{(M^i)}$. 
However, any $\Sigma$ can be dissected in different
ways. The over-used term {\it duality} means here for the
 invertible matrices 
relating the chiral blocks of different dissections. For example,
the left dissection in Figure 2 of the $(g,n)=(0,4)$ surface
corresponds to a matrix 
$F\left[\matrix{M&N\cr L&P}\right]$ of size $n\times n$ for $n=
{\rm dim}\,\frak{F}^{(0,4)}_{(L,M,N,P)}$ (for an appropriate
orientation of surfaces and punctures -- a minor technicality we've been 
ignoring), called the fusing matrix.
Likewise, the right dissection defines the braiding 
matrix $B$.

\medskip\epsfysize=.8in \centerline{\epsffile{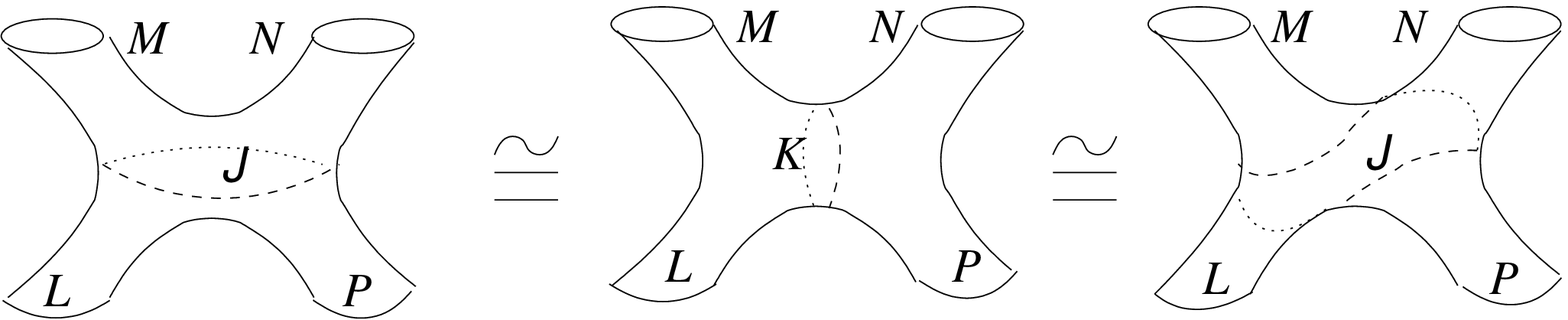}}\medskip

\centerline{{\smcap Figure 2.} The fusing and braiding matrices $F,B$}\medskip

All duality transformations are  built up from a few elementary
 ones, like $B$ and $F$. By decomposing surfaces
in different ways, we get relations between these elementary dualities.
For example, the $B$-matrices obey an equation of the form $BBB=BBB$ called
the Yang--Baxter equation, and this is the source of its name `braiding'.

Consider  four marked points $w_i$ on the sphere $\C\cup\{\infty\}$.
Using the M\"obius(=conformal) symmetry of the sphere,
move $w_i$ to $0,w,1,\infty$, respectively, where $w$ is the cross-ratio
${(w_1-w_2)(w_3-w_4)\over (w_1-w_3)(w_2-w_4)}$.
If we
label all four marked points 
with the Ising module $\sigma\in\Phi$, 
then the space ${\frak F}^{(0,4)}_{(\sigma\sigma
\sigma\sigma)}$ of chiral blocks is two-dimensional; 
choosing the smallest energy state in $\sigma$, the chiral blocks are spanned by
$$\eqalignno{{\cal F}_1(w)=&\,{\sqrt{1+\sqrt{1-w}}\over (w(1-w))^{1/8}} 
\ ,&(4.4)\cr
{\cal F}_2(w)=&\,{\sqrt{1-\sqrt{1-w}}\over (w(1-w))^{1/8}}
\ .&}$$
Fusing interchanges $w_1=0$ and $w_3=1$, hence involves the M\"obius transformation
$w\mapsto (1-w)/(1-0)=1-w$. Likewise, braiding interchanges $w_2$ with $w_3$,
and so involves  $w\mapsto (0-1)/(0-w)=1/w$.
The braiding and fusing matrices here become
$$\eqalignno{B\left[\matrix{\sigma&\sigma\cr\sigma&\sigma}\right]&\,
={e^{-\pi \i/8}\over \sqrt{2}}\left(\matrix{1&\i\cr \i&1}\right)\ ,&(4.5)\cr
F\left[\matrix{\sigma&\sigma\cr\sigma&\sigma}\right]&\,
={1\over \sqrt{2}}\left(\matrix{1&1\cr 1&-1}\right)\  .&}$$

\noindent{{\bf 5. Monodromy in RCFT.}} 

The fractional powers in the Ising blocks (4.4) tell us they
have branch-point singularities -- we must make cuts in the $w$-plane
to get holomorphic functions there. If instead 
we analytically continue these functions along a closed curve, the value
of the block need not return to the same value. For example,
 consider the circle
$w(t)=r\,e^{2\pi\i t}$ for $r$ small: the value of 
${\cal F}_i(w)$ at $t=1$ is $e^{-\pi\i/4}$ times its value at $t=0$. 
This factor $e^{-\pi\i/4}I$ is the {monodromy} about $w=0$. Likewise,
 their monodromy about $w=1$ is 
$$\left(\matrix{{\cal F}_1(w)\cr{\cal F}_2(w)}\right)\mapsto\left(\matrix{
0&e^{-2\pi\i/8}\cr e^{-2\pi\i/8}&0}\right)\left(\matrix{{\cal F}_1(w)\cr
{\cal F}_2(w)}\right)\ .$$

Reintroducing the four 
coordinates $w_i$, the chiral blocks ${\cal F}_i$ will be
holomorphic on  the universal cover
 of the configuration space $\frak{C}_4(S^2)=\{(w_1,w_2,w_3,w_4)\in(S^2)^4\,|\,
 w_i\ne w_j\}$ of the Riemann sphere. Analytically
continuing along any closed path $\gamma$ in $\frak{C}_4(S^2)$ 
defines an action of the fundamental group $\pi_1(\frak{C}_4(S^2))$ 
-- the pure braid group of the sphere with four strands --
 on the space ${\frak F}^{(0,4)}_{(\sigma\sigma\sigma\sigma)}$ of chiral blocks.
 For example, the monodromy about $w=1$ found above
corresponds to the pure braid $\sigma_{1w}^2$, where $\sigma_{1w}$ is the
twist of the 1-strand with the $w$-strand (each strand corresponds to one
of the points $0,1,w,\infty$).  Actually,  the full spherical braid group 
${\frak B}_4(S^2)$ acts:  $\beta\in{\frak B}_4(S^2)$ maps the space
${\frak F}_{(M^1,M^2,M^3,M^4)}^{(0,4)}$ to ${\frak F}^{(0,4)}_{(M^{\beta 1},
M^{\beta 2},M^{\beta 3},M^{\beta 4})}$, where $\beta i$
is the associated permutation. For example, the twist $\sigma_{1w}\in{\frak B}_4(S^2)$
is the braiding matrix $B\left[\matrix{\sigma&\sigma\cr\sigma&\sigma}\right]$.

Equivalently, as a
`function' on the configuration space, the chiral blocks form  
holomorphic horizontal sections of a projectively flat vector bundle. What this means is 
that each chiral block satisfies a system of partial differential equations 
(the {\it Knizhnik--Zamolodchikov} or 
KZ equations) describing how to parallel-transport it around
configuration space, and flatness says it {\it locally} depends only
on the moduli space parameters (and not on the path chosen). {\it Globally},
however, there will be monodromy.

More generally, a chiral block ${\cal F}$ for an enhanced surface $\Sigma$ 
is a multi-valued function on the appropriate moduli space.
To make it well-defined, 
${\cal F}$ can be lifted to the corresponding Teichm\"uller space. There
will be an action of the corresponding mapping class group $\widehat{\Gamma}_{g,n}$, coming from
monodromy. In other words, the space ${\frak F}_{(M^i)}^{(g,n)}$ of chiral 
blocks carries a representation $\rho^{(g,n)}_{(M^i)}$ of $\widehat{\Gamma}_{g,n}$.
This $\widehat{\Gamma}_{g,n}$-representation is built up from the duality 
matrices, such as the braiding and fusing matrices. As we shall see below,
this picture unifies the Jones knot polynomial, the modularity of Monstrous
Moonshine, and many other phenomena.

For example, we can dissect the torus using a single vertical cross-sectional
cut, or using a
horizontal equatorial cut; one basis is $\chi_M(\tau)$ and the other is
$\chi_M(-1/\tau)$. Duality says  that they both span the same space:
$$\chi_M(-1/\tau)=\sum_{N\in\Phi} S_{MN}\,\chi_N(\tau)\ .\eqno(5.1a)$$
Likewise, performing a Dehn twist about the vertical cut, we obtain
$$\chi_M(\tau+1)=\sum_{N\in\Phi} T_{MN}\,\chi_N(\tau)\ .\eqno(5.1a)$$
Here, $S,T$ are complex matrices. Together, $\tau\mapsto-1/\tau$ and
$\tau\mapsto\tau+1$ generate the modular group PSL$_2(\Z)$, and
$S,T$ generate a true representation $\rho^{(1,0)}$ of the central extension 
SL$_2(\Z)=\Gamma_{1,0}$. 
Hence RCFT characters $\chi_M$ form a weight-0 vector-valued modular form for
SL$_2(\Z)$, with multiplier $\rho^{(1,0)}$.

For example, the matrix $T$ for WZW models
involves the quadratic Casimir of $G$, while the matrix $S$ 
involves characters of $G$ evaluated at elements of finite order. 
For the Ising model, these matrices are
$$S={1\over 2}\left(\matrix{1&1&\sqrt{2}\cr 1&1&-\sqrt{2}\cr 
\sqrt{2}&-\sqrt{2}&0}\right)\, ,\
T=\left(\matrix{e^{-\pi\i/24}&0&0\cr 0&-e^{-\pi\i/24}&0\cr 0&0&e^{\pi\i/12}}
\right)\eqno(5.2) $$

Perhaps the most elegant treatment of the finite-dimensional representations
of a compact Lie group $G$ is  Borel--Weil theory, which constructs them
via the $G$-action on line bundles over the flag manifold
$G_\C/B$. Something similar
happens to the Virasoro algebra $\frak{Vir}$, with now the moduli spaces of 
curves playing the role of the flag manifold and mapping class groups taking 
the place of the Weyl group. A copy
of $\frak{Witt}=Vect(S^1)$ attached to the $i$th puncture on an enhanced
surface of type $(g,n)$ acts naturally on the moduli 
space $\widehat{{\frak M}}_{g,n}$:
the vector field $z_i^\ell\partial/\partial z_i$ for $\ell\ge 1$ changes 
the local coordinate $z_i$; $\partial/\partial z_i$ moves the puncture; and
 $z_i^\ell\partial/\partial z_i$ for $\ell\le -1$ can change the
conformal structure of the surface. This infinitesimal action fills out the 
tangent space to any point on $\widehat{{\frak M}}_{g,n}$. In this picture,
the central extension of $\frak{Witt}$ to $\frak{Vir}$ arises geometrically 
as a curvature effect. The KZ equations say  roughly that the desired
sections of $\widehat{{\frak M}}_{g,n}$-vector bundles should  respect this 
$\frak{Vir}$ action.

\bigskip\noindent{\bf 6. Correlation functions.}

Our primary interest is the ${{\Gamma}}_{g,n}$-action $\rho^{(g,n)}_{(M^i)}$
on the spaces
${\frak F}_{(M^i)}^{(g,n)}$. Correlation functions -- the quantities of
physical interest -- are sesquilinear combinations of chiral blocks which have
trivial monodromy. In a typical quantum field theory the correlation functions
are computed perturbatively, but in RCFT they can be found exactly.

For example, the toroidal {correlation function} is
$${\cal Z}(\tau):=\sum_{M\in\Phi,
\overline{N}\in\overline{\Phi}} {\cal Z}_{M,\overline{N}}\,
\chi_M(\tau)\,\chi_{\overline{N}}(\overline{\tau})\ \eqno(6.1)$$
(recall (3.4)). It is required to be invariant under the SL$_2(\Z)$ action
(5.1) on the chiral blocks $\chi_M(\tau)$ of the torus. For the special
case of the Ising model (recall (5.2)), the unique solution to the various
constraints is
$${\cal Z}(\tau)=\chi_{{\cal V}}(\tau)\,{\chi_{{\cal V}}(\overline{\tau})}+
\chi_\epsilon(\tau)\,{\chi_\epsilon(\overline{\tau})}+
\chi_\sigma(\tau)\,{\chi_\sigma(\overline{\tau})}\ .\eqno(6.2)$$
Its SL$_2(\Z)$-invariance follows from the unitarity of
the matrices (5.2). The classification of possible toroidal correlation
functions for WZW models involves quite interesting Lie theory and 
number theory -- see e.g.  [17] for a short proof of 
Cappelli--Itzykson--Zuber's A-D-E classification of SU$_2(\C)$. 

The most elegant and general construction of correlation functions from
chiral blocks uses topological
field theories and the language of category theory [14]. 

\bigskip\noindent{\bf 7. Genus 0: braids and knots.}

In the 1980s, Jones studied the combinatorial characterisation of embedding
one von Neumann algebra (a factor) in another, and as an unexpected by-product
obtained new representations of braid groups ${\cal B}_n$. We can obtain a
knotted link from a braid by gluing the $n$ top endpoints of the braid to
the corresponding bottom ones. It is possible to characterise completely
(using the `Markov moves') the different braids which yield the same knot,
and remarkably Jones' ${\cal B}_n$-representations respect this redundancy,
in the sense that Jones could obtain from his representations (using the
trace in the underlying von Neumann algebra), a polynomial knot invariant [22].

Witten [32] reinterpreted Jones' braid group representations as due to the 
(projective) representation of the genus-0 mapping class groups $\Gamma_{0,n}$, 
coming from the SU$_2(\C)$ WZW model (as we shall see, $\Gamma_{0,n}$ is 
essentially a braid group). Witten showed how the $\Gamma_{0,n}$-representation
of any other RCFT similarly gives rise
to other knot invariants, thus generalising Jones' invariant considerably by embedding it
naturally in a much broader context. von Neumann algebras arise in quantum field theory
(hence RCFT) through the assignment to each region of space-time
of the observables measurable in that region;
when one region is a subset of another, then its algebra of
observables is embedded in the other.

$\Gamma_{0,n}$ has presentation
$$\eqalignno{\Gamma_{0,n}=\langle \sigma_1,\ldots,\sigma_{n-1}\,|&\,
\sigma_i\sigma_j=\sigma_j\sigma_i\
(|i-j|>1),\ \sigma_i\sigma_{i+1}\sigma_i=\sigma_{i+1}\sigma_i\sigma_{i+1},&\cr
&\sigma_1\cdots
\sigma_{n-1}\sigma_{n-1}\cdots\sigma_1=1=(\sigma_1\cdots\sigma_{n-1})^4\rangle
&\cr&\cong{\cal B}_n
(S^2)/\Z_2\ ,&(7.1)}$$
the quotient of the spherical braid group by its centre. RCFT obtains the 
$\Gamma_{0,n}$-representation $\rho^{(0,n)}_{(M^i)}$
by assigning each generator $\sigma_i$ to a matrix in block-diagonal form, whose
blocks are braiding matrices. We get a different representation, of dimension given by
(4.1), for every choice $M^1,\ldots,M^n$ of ${\cal V}$-modules. An element $\beta\in\Gamma_{0,n}$
sends the space ${\frak F}^{(0,n)}_{(M^i)}$ to ${\frak F}^{(0,n)}_{(M^{\beta i})}$,
so to get a representation of the full $\Gamma_{0,n}$ we should sum
${\frak F}^{(0,n)}_{(M^i)}$ over all reorderings of $(M^i)$.
It is common to lift this $\Gamma_{0,n}$-action to the
braid group ${\cal B}_n$ in the manner clear from (7.1). In all known examples
it seems, these representations are always defined over some cyclotomic field.

The groups $\Gamma_{0,0}\cong\Gamma_{0,1}\cong 1$, $\Gamma_{0,2}\cong\Z_2$, $\Gamma_{0,3}
\cong S_3$ are all finite and so aren't very interesting. This is because the M\"obius
symmetry on $S^2$  is triply transitive. However, the image of $\Gamma_{0,n}$ for $n\ge 4$
will usually be infinite -- e.g. in the special case of two-dimensions,
the question of $\Gamma_{0,4}$ having finite image reduces to Schwarz' 
classical analysis of the finite monodromy of the hypergeometric equation,
and as such is very rare.

$\Gamma_{0,4}$ is an extension of PSL$_2(\Z)$ by $\Z_2\times\Z_2$, and the part of 
$\Gamma_{0,4}$ corresponding to a trivial permutation $\beta i=i$ is isomorphic
 to the principal congruence subgroup $\Gamma(2)/
\pm 1$. Using the $\Gamma(2)$-Hauptmodul $\theta_2(\tau)^4/\theta_3(\tau)^4$, we can lift
the chiral blocks ${\cal F}(w)$ to the upper-half plane $\H$, and in this way interpret these chiral
blocks as vector-valued modular forms for SL$_2(\Z)$. For example, the
Ising blocks (4.4) would now become
$$\eqalignno{{\cal F}_1(\tau)=&\,q^{-{1\over 16}}(1+q^{{1\over 2}}+3q+4q^{{3\over 2}}+5q^2
+8q^{{5\over 2}}+11q^3+\cdots)\ ,&\cr
{\cal F}_2(\tau)=&\,q^{{3\over 16}}(2+2q^{{1\over 2}}+2q+4q^{{3\over 2}}+8q^2+10q^{{5\over 2}}+12q^3+\cdots)\ .&(7.2)}$$
These lifts ${\cal F}(\tau)$ will always be holomorphic in $\H$, but can have poles
at the cusps. The weight can be rational because the $\Gamma(2)$-representation will
typically be projective; in fact arbitrary rational weight is possible. The theory of these
vector-valued modular forms of arbitrary rational weight for PSL$_2(\Z)$, and no restriction
on the kernel of the multiplier, has been developed recently [4], [26] and is quite rich.
It is interesting that RCFT produces plenty of examples of these.

The $\Gamma_{0,n}$-action coming from WZW models is especially interesting.
Consider SU$_2(\C)$ for concreteness. 
Choose $n$ distinct points $z_1,\ldots,z_n\in\C$ and $n$ $\frak{sl}_2(\C)$-modules
$V^i$; write $\widehat{V}^i$ for the corresponding $\frak{sl}_2^{(1)}$-modules. Then 
the conformal blocks ${\cal F}\in {\frak{F}}^{(0,n)}_{(\hat{V}^i)}$ are
precisely the functions ${\cal F}:\frak{C}_n(S^2)\rightarrow V_1\otimes\cdots\otimes 
V_n$ satisfying the {\it KZ equations} [27] 
$${\partial {\cal F}\over\partial z_i}={1\over k+2}\sum_{j\ne i}{\Omega_{ij}\over 
z_i-z_j}\,{\cal F}\ ,\eqno(7.3)$$
where 
$\Omega_{ij}/(z_i-z_j)$ is the classical Yang--Baxter $r$-matrix for SU$_2(\C)$.
As mentioned earlier, any solution to (7.3) can be parallel-transported
through $\frak{C}_n(S^2)$; projective flatness means that this
 parallel-transport along a closed loop depends (up to a projective
factor) only on the homotopy-class of the loop. In other words,
the space ${\frak{F}}^{(0,n)}_{(\hat{V}^i)}$ of solutions to (7.3) carries a
projective reprsentation of the pure spherical braid group
$\pi_1(\frak{C}_n(S^2))$. The Drinfel'd--Kohno monodromy theorem expresses
this monodromy in terms of the $6j$-symbols of the quantum group $U_q(\frak{sl}_2(\C))$,
for $q=e^{\pi \i/(k+2)}$, which are straightforward to compute [25]. Something
similar happens for any $G$.

The infinitely many irreducible finite-dimensional modules of a simple Lie 
algebra $\frak{g}$ naturally span a symmetric monoidal category
(see [31] for definitions); its representation ring is isomorphic to
a polynomial ring in $r$ variables, where $r$ is the rank of the algebra. On 
the other hand, the finitely many level $k$ irreducible integrable 
modules of the affine algebra $\frak{g}^{(1)}$ span (among other things) a
braided monoidal category; the corresponding representation 
ring is called a {\it fusion ring}{} and has structure constants equal to
the fusion coefficients (3.6). 
The key ingredient in this category -- the braiding -- comes from the
braid group monodromy of (7.3). Something similar happens for any RCFT.

\bigskip\noindent{\bf 8. Genus 1: modularity.}

The characters $\chi_\la$ of the affine algebra $\frak{g}^{(1)}$ are defined 
exactly as for semi-simple $\frak{g}$, as
a sum of exponentials of the  Cartan
subalgebra, though the sum will now be infinite.
In fact a miracle happens: the character $\chi_\la$ will be a
modular function for some subgroup of SL$_2(\Z)$! One of the coordinates
of the Cartan subalgebra of ${\frak{g}}^{(1)}$ plays the role of $\tau\in\H$,
and the others come along for the ride. 
The algebraic proof of this modularity makes it
look accidental: the character $\chi_\la$ is expressed as a fraction;
the denominator is automatically a modular form for SL$_2(\Z)$, by the simple
combinatorics of affine algebras;  the numerator is a modular form (in fact a
lattice theta function) for some congruence
subgroup of SL$_2(\Z)$, because the Weyl group of ${\frak{g}}^{(1)}$
contains translations in a lattice; their quotient yields a modular function.

RCFT provides a much more satisfying explanation for this
unexpected modularity.  The characters $\chi_\la(\tau)$ of these
affine algebra modules equal the chiral blocks (4.2) of the corresponding
WZW model, and the action (5.1) of $\Gamma_{1,0}={\rm SL}_2(\Z)$ coming from 
RCFT explains their unexpected 
modularity. This relation with RCFT also tells us the 
${{\frak g}}^{(1)}$-modules are simultaneously ${\frak{Vir}}$-modules.
All of this was known to algebraists before the
relation of affine algebras to RCFT was developed, but this relation emphasises
that these properties of affine algebra modules are not accidental
but naturally fit into a much broader perspective.

There is more to being a modular form or function 
than transforming nicely with respect to SL$_2(\Z)$. Good behaviour at the
cusps of $\H$ is also crucial, as they compactify the domain. These cusps
correspond to a pinched torus; their analogue
for the other moduli spaces are surfaces with nodes (this is the Deligne--Mumford
compactification). RCFT requires nice behaviour (`factorisation') 
of chiral blocks as we move in moduli space toward these degenerate
surfaces. This connects the moduli spaces of different topologies, and
tells us RCFT is naturally defined on a `universal tower' of moduli spaces.

Interesting modularity certainly isn't restricted to affine algebras. Indeed,
the Monstrous Moonshine conjectures (see e.g. [18]) relate character
values of the monster finite simple group to various Hauptmoduls. For example,
the first nontrivial coefficient (196884) in the $j$-function  nearly
equals the dimension 196883 of the first nontrivial representation of the
monster. Conjecturally, Hauptmoduls are associated to pairs of commuting
elements in the Monster -- e.g. the $j$-function is assigned to $(e,e)$.
The starting point to our (still 
incomplete) understanding of these conjectures is the construction [13] 
of an RCFT  (with central charge $c=24$, 
$\|\Phi\|=1$, and anti-holomorphic
chiral algebra $\overline{{\cal V}}=\C$) whose symmetry group equals the 
Monster and
whose single character  (4.2) equals the $j$-function. 

To some crude extent, Moonshine can then be interpreted as the conjunction of
two different pictures of quantum field theory, applied to that very special 
RCFT: the Hamiltonian picture, which provides us a
Hilbert space (state-space) carrying an action of the Monster, and an energy
operator $L_0$ such that (4.2) is defined; and the Feynman picture, which
lives in moduli space and which makes modularity manifest.
As explained at the end of Section 5,  the Virasoro algebra, through its 
action on the moduli spaces $\widehat{\frak{M}}_{g,n}$,
 lies at the heart of Moonshine.

Can we see more directly why the
RCFT characters $\chi_M(\tau)$ of (4.2) should have anything to do with 
modularity? The chiral blocks on the torus
can be obtained from those of the plane $\C$, by first
considering the map $z\mapsto w:=e^{2\pi\i z}$. Though holomorphic,
it changes the global topology, 
sending the plane $\C$ to the annulus $\C\backslash\{0\}$, and this topology
change is responsible for the $-c/24$ in (4.2). To
obtain our torus, we now identify $w$ and $qw$, where as always $q=e^{2\pi\i\tau}$.
This is equivalent to taking the finite annulus $\{w\in\C\,|\,|q|<|w|<1\}$
and sewing together its two boundary circles appropriately.
The resulting torus is conformally equivalent to $\C/(\Z+\Z\tau)$.
Applying this construction to chiral blocks, we find that those 
for the torus are indeed given by (4.2) (e.g. the trace comes from sewing).
The proof [33] of modularity of vertex operator algebra characters follows
this outline.

The SL$_2(\Z)$-representation (5.1) is defined over a cyclotomic field for any 
RCFT [8], and its kernel contains a congruence subgroup [3]. Perhaps the
latter isn't so surprising, considering that ${\frak{F}}^{(1,0)}$ has a basis 
(4.2) with integer $q$-expansions. Intimately connected with this 
congruence subgroup property, the matrices $S,T$ have nice properties with 
respect to the cyclotomic Galois group [8,3].

$\Gamma_{1,1}$ is also SL$_2(\Z)$, and so the chiral blocks in 
${\frak{F}}^{(1,1)}_M$ also form vector-valued modular forms for SL$_2(\Z)$.
The weight though can be arbitrary rational numbers and the kernel need
not be of finite index. So together with the ${\cal F}\in\frak{F}^{(0,4)}_{(
M^i)}$ (as explained last section), RCFT is a rich source of vector-valued
modular forms.
These $\Gamma_{1,1}$-representations are explicitly known in terms of the
duality matrices (see e.g. [30]), and so the machinery of [4] allows these
chiral blocks to be explicitly found [5].

But specialising to the WZW models, we should expect a nice Lie theoretic 
answer, and indeed the complete answer is known for SU$_2(\C)$ [24,11].
The matrix representing $\tau\mapsto\tau+1$ will again be given by quadratic 
Casimirs, but the matrix representing $\tau\mapsto -1/\tau$ involves the
continuous $q$-ultraspherical polynomials, also
 known as the Macdonald polynomials for $\frak{sl}_2(\C)$ (and as such are a
natural generalisation of SU$_2(\C)$-characters, which describe $\tau\mapsto
-1/\tau$ for $\rho^{(1,0)}$). Similarly, the chiral blocks in 
${\frak{F}}^{(1,1)}_\lambda$ can be expressed using Macdonald `polynomials'
for $\frak{sl}_2^{(1)}$ (a natural generalisation of the 
$\frak{sl}_2^{(1)}$-characters which are the chiral blocks in 
$\frak{F}^{(1,0)}$). Something similar can be expected for the other 
WZW models on the punctured torus.
Likewise, evaluating these chiral blocks for the Moonshine RCFT should
extend Moonshine to noncommuting pairs $g,h$ in the monster [5].

The space $\frak{F}_{{\cal V}}^{(1,1)}$, with the ${\cal V}$-module ${\cal V}$,
 can be identified with $\frak{F}^{(1,0)}$
except that in the former we have the freedom to evaluate the block at
any state $v\in{\cal V}$. In particular, taking $v$ to have the minimum energy 
-- the vacuum -- recovers the characters (4.2), but taking other $v$ will
give a vector-modular form of even weight, with the same SL$_2(\Z)$-multiplier 
$\rho^{(1,0)}$ as in (5.1). [10] evaluated these for the Moonshine RCFT and
found that all modular forms (of the right shape) arise.
What is interesting is that the coefficients of these modular forms will
have Moonshine-like interpretations involving characters of the stabiliser
of the state $v$.

\bigskip\noindent{\bf 9. Higher genus.}

For a fixed RCFT, chiral blocks ${\cal F}\in\frak{F}^{(g,n)}_{(M^1,\ldots,M^n)}$
yield vector-valued automorphic functions for the infinite
discrete groups $\Gamma_{g,n}$, each realising a finite-dimensional
representation $\rho^{(g,n)}_{(M^i)}$ of $\Gamma_{g,n}$. This tower of automorphic functions
is coherent in the sense that it respects basic operations like sewing or
pinching the surfaces. As mentioned earlier, $(g,n)=(1,0),(1,1),(0,4)$ all
give a vector-valued modular form for SL$_2(\Z)$; for $(1,0)$ this is  
 a classical object, being weight-0 and invariant under some 
congruence subgroup, but for $(1,1)$ and $(0,4)$ the weight is rational
and the image of the multiplier will usually be infinite.
Relatively little is known in
higher genus $g$, and surely it is a direction for important future
research. The main open challenge is to identify the special features and
structures occurring here. In this sense most of the work done has been negative.
In this section we sample a few of the highlights.

Most of the work has focused on the kernel and image of these representations
$\rho^{(g,n)}_{(M^i)}$. The groups $\Gamma_{0,n}$, $n\le 3$, are finite;
all other  ker$\,\rho^{(g,n)}_{(M^i)}$ will be infinite, since 
$\Gamma_{g,n}$ is generated by 
infinite-order Dehn twists but $\rho^{(g,n)}_{(M^i)}$ maps each of these
to a finite-order matrix. However, for fixed $n$, the intersection over
all $k$ of the kernel of $\rho^{(g,0)}$ for the SU$_n(\C)$ WZW model
at level $k$, is trivial in any genus $g>2$ [1]. 

The image of $\rho^{(1,0)}$ is always finite [3], but this is atypical:
it is expected that a generic RCFT will have
all other images infinite. For example, Funar [15] found that all
im$\,\rho^{(g,0)}$ will be infinite for SU$_2(\C)$ WZW models at all levels
$k>8$, and all genus $g>1$. Moreover, Masbaum [29] found an infinite-order matrix in 
im$\,\rho^{(0,4)}$ for those RCFT. 

In the RCFT associated to even self-dual lattices $L$ (where the strings live
on the torus $\R^n/L$ for $n={\rm dim}\,L$), the conformal blocks in
${\frak{F}}^{(g,0)}$ can be expressed in terms of Siegel theta functions,
and the Torelli subgroup of $\Gamma_{g,0}$ is in the kernel of $\rho^{(g,0)}$.
This is very atypical for RCFT, e.g. it is known to fail for SU$_2(\C)$ WZW at 
most levels. 

On the other hand, these representations for all known RCFT seem to be
always definable over a cyclotomic field. A notion of integrality for
these representations is being developed [20].

A class of RCFT very conducive to study are the so-called holomorphic orbifolds,
associated to the Drinfel'd double of a finite group $G$. In this case,
the chiral blocks in $\frak{F}^{(g,0)}$ are parametrised by Hom$(\pi_1(\Sigma_g),
G)/G$, i.e. group 
homomorphisms $\varphi$ from $\pi_1$ of a genus-$g$ surface $\Sigma_g$ into $G$,
where we identify $\varphi(\sigma)$ and $g^{-1}\varphi(\sigma)g$.
$\Gamma_{g,0}$ acts naturally
on $\pi_1(\Sigma_g)$ and hence $\rho^{(g,0)}$ here becomes a permutation 
representation. This means im$\,\rho^{(g,0)}$ here is manifestly finite.
However, we can see from this explicitly that ker$\,\rho^{(g,0)}$ won't
contain the Torelli generators listed by [6], at least for generic groups
$G$, and so even in this extremely well-behaved theory, $\rho^{(g,0)}$
doesn't factor through to a representation of Siegel's modular group 
Sp$_{2g}(\Z)$. [12] show
that im$\,\rho^{(0,n)}$ will always be finite here, and it is tempting
to guess that all im$\,\rho^{(g,n)}$ is finite here.

\bigskip\noindent{{\bf Acknowledgements.}} I've benefitted greatly from 
communications with Peter B\'antay, Gregor Masbaum, and Mark Walton. This 
research is supported in part by NSERC.

\bigskip\noindent{\bf Bibliography.}\medskip

\item{[1]} J. E. Andersen, `Asymptotic faithfulness of the quantum $SU(n)$
representations of th mapping class groups', {\it Ann.\ Math.} {\bf 163}
(2006) 347--368.

\item{[2]} {B.\ Bakalov} and {A.\  Kirillov, Jr.}, {\it Lectures on Tensor 
Categories and Modular Functors} (American Mathematical Society, Providence 2001).

\item{[3]} P. B\'{a}ntay, `The kernel of the modular representation 
and the Galois action in RCFT', {\it Commun.\ Math.\ Phys.} {\bf 233} (2003) 
423--438.

\item{[4]} P. B\'{a}ntay and T. Gannon, `Vector-valued modular functions
for the modular group and the hypergeometric equation', arXiv:math/0705.2467.

\item{[5]} P. B\'{a}ntay and T. Gannon, in preparation.

\item{[6]} J. S. Birman, `On Siegel's modular group', {\it Math. Ann.}
{\bf  191} (1971)  59--68.

\item{[7]} {J.\ S.\ Birman}, {\it Braids, Links, and Mapping Class Groups}
(Princeton University Press, Princeton 1974).

\item{[8]} {A.\ Coste} and {T.\ Gannon}, `Remarks on
Galois in rational conformal field theories', {\it Phys.\ Lett.}\ {\bf B323}
(1994) 316--321.

\item{[9]} {P.\ Di Francesco, P.\ Mathieu}, and {D.\ S\'en\'echal}, {\it
Conformal Field Theory} (Springer, New York 1996).

\item{[10]} C. Dong and G. Mason, `Monstrous moonshine at higher weight',
{\it Acta Math.} {\bf 185} (2000) 101--121.

\item{[11]} P. I. Etingof and A. A. Kirillov Jr, `On the affine analogue of
Jack and Macdonald polynomials', {\it Duke Math. J.} {\bf 78} (1995) 229--256.

\item{[12]} P. Etingof, E. Rowell and S. Witherspoon, `Braid group reprsentations
from twisted quantum doubles of finite groups', math/070327.

\item{[13]} {I.\ Frenkel}, {J.\ Lepowsky}, and {A.\ Meurman}, {\it
Vertex Operator Algebras and the Monster} (Academic Press, San Diego
1988).

\item{[14]} {J.\ Fuchs, I.\ Runkel,} and {C.\ Schweigert},
`Boundaries, defects and Frobenius algebras', {\it Fortsch.\ Phys.} {\bf 51}
(2003) 850--855.

\item{[15]} L. Funar, `On the TQFT representations of the mapping class groups',
{\it Pacif. J. Math.} {\bf 188} (1999) 251--274.

\item{[16]} {M.\ R.\ Gaberdiel}, `Introduction to conformal field
theory', {\it Rep.\ Prog.\ Phys.} {\bf 63} (2000) 607--667.

\item{[17]} {T.\ Gannon}, `The Cappelli--Itzykson--Zuber A-D-E classification',
{\it Rev.\ Math.\ Phys.} {\bf 12} (2000) 739--748.

\item{[18]} T.\ Gannon, {\it Moonshine Beyond the Monster} (Cambridge University
Press, 2006).

\item{[19]} {K.\ Gawedzki}, `Conformal field theory', S\'eminaire
Bourbaki, {\it Ast\'erisque} {\bf 177-178} (1989) 95--126.

\item{[20]} P. M. Gilmer and G. Masbaum, `Integral lattices in TQFT', arXiv:
math.QA/0411029.

\item{[21]} {Y.-Z.\ Huang}, {\it Two-Dimensional Conformal Geometry and
Vertex Operator Algebras} (Birkh\"auser, Boston 1997).

\item{[22]} V. Jones and V.\ S.\ Sunder, {\it Introduction to Subfactors}
(Cambridge University Press, 1997).

\item{[23]} {V.\ G.\ Kac}, {\it Infinite Dimensional Lie Algebras}, 
3rd edn (Cambridge University Press, Cambridge 1990).

\item{[24]} A. A. Kirillov Jr, `On an inner product in modular tensor categories',
{\it J. Amer. Math. Soc.} {\bf 9} (1996) 1135--1169.

\item{[25]} A.\ N.\ Kirillov and N.\ Yu.\ Reshetikhin, `Representations of the
algebra $U_q(sl(2))$, $q$-orthogonal polynomials and invariants of links',
in: {\it Infinite-dimensional Lie algebras and groups (Luminy-Marseille, 1988)}
(World Scientific, Teaneck NJ 1989) 285--339.

\item{[26]} M. Knopp and G. Mason, `Vector-valued modular
forms and Poincar\'{e} series', {\it Illinois J. Math.} {\bf 48} (2004) 1345-1366.

\item{[27]} T. Kohno, {\it Conformal Field Theory and Topology} (American Math.
Soc., 2002).

\item{[28]} {J.\ Lepowsky} and {H.\ Li}, {\it Introduction to
Vertex Operator Algebras and Their Representations} (Birkh\"auser, Boston
2004).

\item{[29]} G. Masbaum, `An element of infinite order in TQFT-representations 
of mapping class groups', {\it Contemp. Math.} {\bf 233} (1999) 137--139.

\item{[30]} {G.\ Moore} and {N.\ Seiberg}, `Classical and
quantum conformal field theory', {\it Commun.\ Math.\ Phys.} {\bf 123}
(1989) 177--254.

\item{[31]} V. G. Turaev, {\it Quantum Invariants of Knots and 3-manifolds}
(de Gruyter, Berlin 1994).

\item{[32]} E.\ Witten, `Quantum field theory and the Jones polynomial',
{\it Commun.\ Math.\ Phys.} {\bf 121} (1989) 351--399.

\item{[33]} {Y.\ Zhu}, `Modular invariance of characters
of vertex operator algebras', {\it J.\ Amer.\ Math.\ Soc.} {\bf 9}
(1996) 237--302.

\end